%% file: article.tex
\documentclass[USenglish,twocolumn]{article}
\usepackage{amsmath}
\ifx\directlua\undefined\ifx\XeTeXcharclass\undefined
  \else\RequirePackage[no-math]{fontspec}[2017/03/31]\fi 
  \else\RequirePackage[no-math]{fontspec}[2017/03/31]\fi 

\usepackage[big]{dgruyter}
\usepackage{microtype}
\input{commands.tex}

\begin{document}

\articletype{Methods}

 \author*[1]{Maximilian Pierer von Esch}
 \author[2]{Andreas V\"olz}
 \author[2]{Knut Graichen} 
 \runningauthor{Pierer v. Esch}
 \affil[1]{Chair of Automatic Control,
  	Friedrich-Alexander-Universit\"at Erlangen-N\"urnberg,
  	e-mail: maximilian.v.pierer@fau.de. Funded by the German Research
Foundation under project no. 464391622.}
 \affil[2]{Chair of Automatic Control,
  	Friedrich-Alexander-Universit\"at Erlangen-N\"urnberg.}
  \title{An Overview of Sensitivity-Based Distributed Optimization and Model Predictive Control}
  \subtitle{Ein Überblick zur sensitivitätsbasierten verteilten Optimierung und modellprädiktiven Regelung}
  \abstract{This paper presents a concise overview of sensitivity-based methods for solving large-scale optimization problems in distributed fashion. The approach relies on sensitivities and primal decomposition to achieve coordination between the subsystems while requiring only local computations with neighbor-to-neighbor communication. We give a brief historical synopsis of its development and apply it to both static and dynamic optimization problems. Furthermore, a real-time capable distributed model predictive controller is proposed which is experimentally validated on a coupled watertank system.}
  \received{...}
  \accepted{...}
  \journalname{...}
  \journalyear{...}
  \journalvolume{..}
  \journalissue{..}
  \startpage{1}
  \aop
  \DOI{...}

\maketitle

\noindent {\bfseries Kurzfassung: } Dieser Beitrag bietet einen kompakten Überblick zu sensitivitätsbasierten Verfahren für die verteilte Lösung hochdimensionaler Optimierungsprobleme. Das Schema nutzt Sensitivitäten und primale Dekomposition, um die Koordination zwischen den Teilsystemen sicherzustellen, wobei lediglich lokale Berechnungen sowie die Kommunikation mit den jeweiligen Nachbarn erforderlich sind. Zusätzlich zu einer historischen Einordnung wird sowohl die Anwendung auf statische als auch auf dynamische Optimierungsprobleme betrachtet. Darüber hinaus wird ein echtzeitfähiges verteiltes modellprädiktives Regelungsverfahren vorgestellt, welches experimentell an einem gekoppelten Wassertanksystem validiert wird.

\vskip\baselineskip
\noindent {\bfseries Schlagwörter: } verteilte Optimierung, verteilte modellprädiktive Regelung, Dekomposition, Sensitivit\"aten

\section{Introduction} 
Large-scale nonlinear programs (NLP) appear in many modern engineering applications, ranging from statistical or machine learning with large amounts of data and many features \cite{Boyd} to large-scale infrastructure systems such as electrical power grids \cite{Molzahn}. These problems are typically high-dimensional, often constrained, and characterized by structured couplings between subsystems or features. Their size and complexity make centralized solution approaches computationally demanding or even infeasible, which has motivated the development of distributed optimization methods \cite{Bertsekas}. By decomposing the central problem into smaller, parallelizable subproblems, such methods enable local computation combined with limited communication, offering scalability, flexibility, and robustness advantages over centralized approaches.

A particularly challenging application domain of large-scale optimization is nonlinear distributed model predictive control (DMPC) \cite{Christofides} or distributed moving horizon estimation (DMHE) \cite{Farina}. Here, each subsystem solves a local optimal control problem (OCP) subject to dynamics and constraints, while considering neighboring subsystems via coupling conditions. The resulting global OCP is nonlinear, constrained, potentially high-dimensional and must be solved within a given sampling time. Distributed optimization techniques are well suited in this case as they enable each subsystem to compute control inputs locally while preferably communicating only with neighboring subsystems. Compared to centralized model predictive control (MPC), this leads to better scalability and robustness. If cooperative approaches are pursued, the same control quality and stability properties as in MPC may be achieved~\cite{Muller}. 

Consequently, a broad variety of approaches with different application areas in mind have been developed. Common methods for constrained convex problems include the alternating direction method of multipliers (ADMM)~\cite{Boyd}, dual decomposition \cite{Everett}, distributed projected gradient descent \cite{Xi} or non-smooth Newton-type schemes \cite{Frasch}. Non-convex problems are significantly more challenging. A straightforward approach for distributed non-convex optimization is to distribute the internal computations of classical nonlinear optimization methods. This procedure has resulted in many bi-level algorithms such as using ADMM to solve the quadratic programs arising in sequential quadratic programming (SQP)~\cite{Stomberg}, distributing the Newton step in interior point methods via ADMM \cite{Engelmann3}, or exploiting the locally convex structure of augmented Lagrangian functions to apply dual decomposition \cite{Bertsekas3}. Another research avenue concerns combining SQP and ADMM resulting in the augmented Lagrangian direction inexact Newton (ALADIN) method~\cite{Houska}.

The focus of this work is on a recently developed sensitivity-based distributed programming (SBDP) method for distributed, non-convex optimization \cite{Pierer, Pierer4}. The approach is fundamentally different from duality-based approaches such as ADMM as it relies on primal decomposition and the incorporation of so-called sensitivities in the local problems. These sensitivities capture the interaction between subsystems and allow for a cooperative solution of the centralized problem with formal convergence guarantees. Its design avoids the overhead of reformulating the original NLP as a generalized consensus problem and instead iterates directly on the coupled variables. As a result, the approach yields subproblems of smaller dimension, reduces communication requirements, and improves scalability. However, it requires careful handling of coupled constraints.
It should be emphasized that we do not present any new theoretical results in this paper, but aim at providing a concise and application-oriented overview of SBDP. 

The contribution of this paper is threefold: First, we present a historical survey on the conceptual development of the SBDP approach. Second, we provide a structured overview of the SBDP method for the solution of graph-structured problems ranging from its application to static optimization over dynamic optimization to DMPC. Third, we present novel experimental results of sensitivity-based DMPC applied to a coupled water tank system which demonstrates its practical performance and potential.

The remainder of the paper is organized as follows. Section~\ref{sec:historical_development} places SBDP into a historical context and highlights related work. Section~\ref{sec:SBDP_static} discusses the application of SBDP to static optimization problems, while Section~\ref{sec:SBDP_dynamic} extends the approach to dynamic optimization problems. Section~\ref{sec:DMPC} introduces a real-time sensitivity-based DMPC scheme. Section~\ref{sec:Experiment} presents both simulation and experimental results, before Section~\ref{sec:Conclusion} concludes the paper.

\section{A brief history of sensitivity- based distributed optimization}
\label{sec:historical_development}
This section provides a brief historical and conceptual overview of sensitivity-based distributed optimization. The goal is to place SBDP into context of classical parametric programming and decomposition methods, identify intermediary ideas ranging from interaction operators to hierarchical modeling approaches, and contrast it with other yet similar decomposition methods. The narrative below is arranged in chronological order and emphasizes the transition from conceptual ideas and related methods to explicit, provably convergent algorithms for distributed optimization. 

\noindent {\bfseries The 1970-80s: Parametric programming and interaction operators.} Sensitivity-based decomposition methods trace back to classical sensitivity analysis in parametric programming \cite{Fiacco} and hierarchical control of multi-level systems \cite{Mesarovic}. Hereby, sensitivity analysis examines how an optimal solution of a parametric NLP changes when problem parameters vary. The behavior of the primal-dual solution of parametric NLPs with respect to parameter perturbations was extensively studied in the 1970s. In particular, it is proven in \cite{Fiacco} under standard regularity conditions, i.e., linear independence constraint qualification (LICQ), second-order sufficient conditions (SOSC), and strict complementarity slackness (SCS), that the primal-dual solution of an NLP varies smoothly w.r.t.\ the parameters. This allows to apply the implicit function theorem to the Karush-Kuhn-Tucker (KKT) system to obtain the local differentiability of the primal-dual solution mapping and the optimal value function. This line of work establishes the legitimacy of using first-order information to approximate how one subsystem responds to changes in variables of other subsystems and is the formal reason why exchanging sensitivities may be used to approximate coupling effects in distributed settings.

In parallel to these developments, hierarchical control concepts and multilevel system theories were developed that formalize interactions among subsystems  \cite{Mesarovic,Mesarovic2,Takahara}. Mesarovic et al. stated that coordinated optimization can be achieved in two different ways: Either modify the model or modify the goal, i.e., the objective function \cite{Mesarovic}. The SBDP approach can be assigned to the latter which are referred to as goal coordination methods.
Specifically, the notion of goal-interaction operators emerged as an early conceptual precursor to the sensitivities later used in the SBDP approach. Mesarovic et al. formulated a two-level hierarchical framework, where each ``infimal'' subsystem has a local objective and receives inputs from a higher level ``supremal'' coordinator \cite{Mesarovic2}. A goal-interaction operator was defined to capture how the aggregate performance of the system depends linearly on each subsystem’s action \cite{Takahara}. Although this work did not produce an explicit algorithm, it formalizes the idea that local goals may be modified to account for interactions.

Furthermore, during the 1970-80s the decomposition methods matured in optimization and process engineering \cite{Lasdon, Lasdon2, Sobieszczanski}. Classic decomposition techniques, such as Benders decomposition, showed how to split large problems by separating complicating variables \cite{Geoffrion}. Hereby, the subproblem feedback to the master is fundamentally a sensitivity information used to form so-called Bender cuts or prices. Benders decomposition shares mathematical and conceptual similarities with the sensitivity-based approaches, i.e., they both use linear approximations of subproblem influence, but they differ in where the linearization is used and in their communication structure. In systems engineering, sensitivities were found to be useful to evaluate choices in design problems of complex, internally coupled systems, predating multidisciplinary design optimization \cite{Sobieszczanski2}.

\noindent {\bfseries The 1990s and 2000s: Alternative decomposition approaches and hierarchical modelling.}
Further theoretical work in the context of parametric programming generalized sensitivity analysis to degenerate and non-regular problems \cite{Ralph,Kyparisis} and investigated the conditions needed to retain the (directional) differentiability of the optimal solution and value function \cite{Jittorntrum}. 

In parallel to these developments, a distinct but related  distributed optimization scheme pursued proximal decomposition via alternating linearization \cite{Kiwiel, Goldfarb}. These methods combine proximal point regularization with block-coordinate linearization, leading to algorithms where each subproblem is solved approximately via a linearized model, and convergence results were established with techniques from proximal minimization. Conceptually, these approaches share with sensitivity-based schemes the idea of iteratively exchanging linearized coupling information between subproblems.

An important aspect in the historical development of sensitivity-based decomposition methods originates in chemical process systems engineering, where the emphasis during the 1990s shifted toward creating structured and hierarchical representations of complex process models \cite{Marquardt,Marquardt2,Stephanopoulos}. The goal at that time was not distributed optimization as considered here, but the establishment of systematic methodologies for modeling and decomposition that would later provide the conceptual basis for interaction-based coordination schemes. Central to this effort was to model the process as modular units (e.g., reactors, separators, heat exchangers) interconnected via material or energy flows which resembles the graph structured approach taken for example in DMPC. 

\noindent {\bfseries The 2010s: Resurgence of interest in DMPC and DMHE.} 
Building on the ideas of goal-interaction operators and structured decomposition, the first practical sensitivity-based distributed optimization schemes were developed in the context of DMPC and DMHE in the 2010s. Scheu et al. proposed a distributed optimization algorithm where each subsystem uses linearized information of the neighboring dynamics and cost around the current iterate to augment its own cost resulting in a cooperative DMPC scheme \cite{Scheu,Scheu2}. This sensitivity-based coordination iteratively updates the primal-dual variables by solving local modified problems in a parallel Jacobi-like fashion. Convergence conditions were established under suitable assumptions for strictly convex, linear-quadratic problems \cite{Scheu}. However, for nonlinear problems formal convergence guarantees and efficient ways to compute the sensitivities were still open questions \cite{Scheu2}. An experimental validation of the sensitivity-driven DMPC scheme followed in \cite{Alvarado}, where the approach displayed competitive performance compared to other methods in a benchmark application to coupled water tanks. 

Subsequently, the idea was also applied to DMHE \cite{Schneider,Schneider2,Schneider3}. These works addressed problems how to treat coupled constraints via an active set strategy \cite{Schneider3}, always provide convergence guarantees for strongly coupled systems by suitable algorithmic modifications \cite{Schneider} and prove stability of the upper layer MHE scheme for a fixed number of underlying optimizer iterations \cite{Schneider,Schneider2}. However, the focus remained on linear systems.

\noindent {\bfseries The 2020s: Expansion to nonlinear, optimal control, and game-theoretic settings.} 
More recently, the sensitivity-based approach has been pushed beyond linear-quadratic settings to general nonlinear and non-convex optimization problems \cite{Pierer, Pierer4}. In particular, a sensitivity-based approach is developed to solve generic NLPs with only local computations and neighbor-to-neighbor communication while providing local convergence guarantees. 

Moreover, the method was extended from static optimization problems to dynamic optimization problems \cite{Huber,Pierer2,Pierer3}. Here, it was shown that the required sensitivities can be computed efficiently and locally via standard optimal control theory. An interesting observation was that if the local subproblems are solved via gradient-based schemes, the sensitivities are obtained for ``free'', in the sense that no additional computations are necessary.  

Furthermore, the idea of enhancing game-theoretic trajectory planning algorithms with sensitivity information has gained traction in the last years. Specifically, a sensitivity-enhanced iterative best response method (SE-IBR) was developed for autonomous racing and intersection scenarios which is deeply rooted in the ideas of the aforementioned works \cite{Spica,Wang,Yuan,Mayer}. In \cite{Spica} it is proven that the converged iterates satisfy a Nash equilibrium. Notably, the practical usefulness of SE-IBR has been confirmed by simulation and experimental results.

In summary, the historical development of sensitivity-based optimization shows a clear progression: from parametric programming and goal-interaction operators, through hierarchical modeling in process engineering, to modern applications in distributed MPC or game-theoretic control. Across these stages, the central idea remained the same, namely to approximate and coordinate subsystem interactions via sensitivity information. The next section formalizes these ideas and provides explicit algorithms for static and dynamic optimization problems. 
\section{Sensitivity-based distributed static optimization}
\label{sec:SBDP_static}
In this section, we discuss the SBDP method for solving static optimization problems in distributed fashion. The scheme relies on constructing local NLPs for each subsystem augmented with the first-order sensitivity information of neighboring subsystems and takes advantage that this information is compactly represented by the directional derivative of the Lagrangian of the neighboring NLPs w.r.t.\ the subsystem's own optimization variables. We discuss different update rules for the primal-dual solution and present the distributed optimization algorithm.  
\subsection{Problem statement}
We consider NLPs which are structured over an undirected, connected graph $\mathcal{G}=(\mathcal{V}, \mathcal{E})$, where the set of nodes $\mathcal{V} = \{1, \dots, M\}$ represents a collection of subsystems. In the following these subsystems are referred to as agents. The edge set $\mathcal{E} \subset \mathcal{V} \times \mathcal{V}$ encodes the interconnection structure of the subsystems. 
The goal of the agents is to solve the central NLP in a cooperative fashion
\begin{subequations}\label{eq:central_NLP}
	\begin{alignat}{2}
		\min_{\vm x_1,\dots,\vm x_M} &\quad  \sum_{i\in\mathcal V} f_i(\vm x_i, \Ni{x}) \label{eq:central_costFunction}\\
		~\st \quad&\quad \vm  g_i( \vm x_i,\Ni{x}) = \vm 0\,,&&\quad \agents \label{eq:central_equality} \\
		&\quad \vm h_i(\vm x_i, \Ni{x})\leq \vm 0 \,,&&\quad \agents \label{eq:central_inequality}\,,
	\end{alignat}
\end{subequations}
where $\vm x_i \in \mathbb{R}^{n_i}$ is the local optimization vector.  The agents $j \in \mathcal{V}$ which are directly coupled with agent $\agents$ are collected in the set ${\mathcal{N}_i := \{j \in \mathcal{V}\,|\,(i,j) \in \mathcal{E}, i \neq j\}}$. These couplings arise via the neighboring optimization variables $\vm x_j \in \mathbb{R}^{n_j}$ through the objective \eqref{eq:central_costFunction} or the (in)equality constraints \eqref{eq:central_equality} -- \eqref{eq:central_inequality} and are summarized as $\Ni{x}:= [\vm x_j]_{\neighs}$.
Each agent minimizes a local objective function $f_i:\mathbb{R}^{n_i}\times\mathbb{R}^{n_{\mathcal{N}_i}} \rightarrow \mathbb{R}$ subject to equality constraints $\vm g_i:\mathbb{R}^{n_i}\times\mathbb{R}^{n_{\mathcal{N}_i}} \rightarrow \mathbb R^{n_{g_i}}$ and inequality constraints $\vm h_i:\mathbb{R}^{n_i}\times\mathbb{R}^{n_{\mathcal{N}_i}}\rightarrow \mathbb R^{n_{h_i}}$ with $ n_{\mathcal{N}_i} := \sum_{\neighs} n_j$. All functions appearing in the central NLP \eqref{eq:central_NLP} are assumed to be sufficiently smooth. 

We define the central Lagrangian of problem~\eqref{eq:central_NLP} as
\begin{equation}
	L(\vm x, \vm \lambda, \vm \mu) := \sum_{\agents}L_i(\vm x_i, \vm \lambda_i, \vm \mu_i, \Ni{x}) \label{eq:central_Lagrangian}
\end{equation}
with the local Lagrangians
\begin{multline}
	L_i(\vm x_i, \vm \lambda_i, \vm \mu_i, \Ni{x}) \\ :=f_i(\vm x_i, \Ni{x}) + \vm \lambda_i\trans \vm  g_i( \vm x_i,\Ni{x}) + \vm \mu_i\trans  \vm h_i(\vm x_i, \Ni{x}) \label{eq:local_Lagrangian}
\end{multline}
for every $\agents$. The quantities $\vm \lambda_i \in \mathbb{R}^{n_{g_i}}$ and $\vm \mu_i \in \mathbb{R}^{n_{h_i}}$ are the Lagrange multipliers associated with the constraints \eqref{eq:central_equality} -- \eqref{eq:central_inequality}, respectively. The centralized decision variable is 
$\vm x = [\vm x_i]_{\agents} \in \mathbb{R}^n$, 
with stacked multipliers 
$\vm \lambda = [\vm \lambda_i]_{\agents} \in \mathbb{R}^{n_g}$ 
and 
$\vm \mu = [\vm \mu_i]_{\agents} \in \mathbb{R}^{n_h}$. 
Together, they form the primal-dual solution vector
$
\vm p := [\vm x^\top,\, \vm \lambda^\top,\, \vm \mu^\top]^\top \in \mathbb{R}^p,
$
of NLP~\eqref{eq:central_NLP} with total dimension $p = n + n_g + n_h$.
For large-scale systems with many decision variables $n\gg n_i$, the centralized NLP~\eqref{eq:central_NLP} quickly becomes intractable and is more efficiently solved in distributed manner. We therefore focus on its solution only via local computations and neighbor-to-neighbor communication.

\subsection{Local problems, sensitivities and primal-dual update}
The structure of NLP \eqref{eq:central_NLP} is taken into account by constructing decoupled, local NLPs in terms of a local search direction $\vm s_i \in \mathbb{R}^{n_i}$ for each $\agents$ which are subsequently solved in each iteration $q=0,1,\dots$ of the SBDP method
	\begin{subequations}\label{eq:local_NLP}
		\begin{alignat}{2}
			\min_{\vm s_i} \quad & \bar f_i\indq(\vm s_i) \label{eq:local_costFunction}
			\\  \st \quad&  \vm{\bar g}_i\indq(\vm s_i)= \vm 0\label{eq:local_equality}  \\
			 & \vm{\bar h}_i \indq(\vm s_i)\leq\vm 0 \label{eq:local_inequality}
		\end{alignat}
	\end{subequations}
with the modified cost functions $\bar f_i: \mathbb{R}^{n_i} \rightarrow \mathbb{R}$, defined as
	\begin{equation} \label{eq:definition_costFunction}
	\bar f_i\indq(\vm s_i) := f_i(\vm x_i\indq + \vm s_i, \Ni{x}\indq) +\frac{\rho}{2}\|\vm s_i\|^2\! +\! \sum_{\neighs}\!\! \nablax\trans L_j\indq \vm s_i
	\end{equation}
	and the local equality and inequality constraints
	\begin{equation}\label{eq:definition_constraints}
		\vm{\bar g}_i\indq(\vm s_i) := \vm g_i( \vm x_i\indq \!+\! \vm s_i,\Ni{x}\indq),\, \vm{\bar h}_i\indq(\vm s_i) := \vm h_i( \vm x_i\indq \!+\! \vm s_i,\Ni{x}\indq)
	\end{equation}
 which depend explicitly on the search direction $\vm s_i$. The implicit dependency of \eqref{eq:definition_costFunction} -- \eqref{eq:definition_constraints} on the neighboring primal and dual variables is captured by the superscript $q$ which indicates that these functions change in each iteration. The objective \eqref{eq:definition_costFunction} consists of three distinct parts. The first term represents the agent's local objective in direction $\vm s_i$. The second is a quadratic regularization term $\frac{\rho}{2}\|\vm s_i\|^2$, where $\rho \geq0$ is a suitable penalty parameter. The third is the sensitivity term which is defined as the directional derivative of neighboring agents' local Lagrangians \eqref{eq:local_Lagrangian} and accounts for the first-order influence of a step in direction $\vm s_i$ on the neighboring objectives. It is given by the gradient of $L_j(\cdot)$ w.r.t.\ $\vm x_i$ evaluated at iteration $q$, i.e., 
 \begin{equation} \label{eq:gradient_abbr}
 \nablax L_j\indq:= \nablax L_j(\vm x_j\indq, \vm \lambda_j\indq, \vm \mu_j\indq, \Nj{x}\indq)\,,
 \end{equation} 
 in direction of $\vm s_i$. 
 Formally, the sensitivity corresponds to the directional derivative of the optimal value function of a neighboring nonlinear program, when $\Nj{x}$ are regarded as parameters. This establishes a direct connection to classical sensitivity analysis in parametric programming~\cite{Fiacco}.
 The required gradient \eqref{eq:gradient_abbr} is computed as 
\begin{align}\label{eq:gradient_withoutneighboraffine}
		& \nablax L_j(\vm x_j, \vm \lambda_j, \vm \mu_j, \Nj{x}) =  \nablax f_j(\vm x_j, \Nj{x}) \nonumber                        \\
		& \phantom{=}+ \nablax\trans \vm g_j(\vm x_j, \Nj{x}) \vm \lambda_j + \nablax\trans \vm h_j(\vm x_j, \Nj{x}) \vm \mu_j\,,
	\end{align}
where $\vm \lambda_j$ and $\vm \mu_j$ are the Lagrange multipliers of neighbor $\neighs$. However, calculating this gradient may require variables from second-order neighbors, which are not accessible in a typical neighbor-to-neighbor communication network. Therefore, each agent $ \agents$ computes the ``mirroring'' gradient $\nablaxj L_i\indq$ and sends it to the respective neighbors $\neighs$. The actual decoupling of \eqref{eq:definition_costFunction} -- \eqref{eq:definition_constraints} is achieved in a primal decomposition fashion by treating the neighboring variables as fixed at the current iteration. 

After solving the local NLP \eqref{eq:local_NLP}, each agent updates its local component $\vm p_i \in \mathbb{R}^{p_i}$, $p_i = n_i + n_{g_i} + n_{h_i}$ of the central primal-dual vector $\vm p$.
To this end, different updates are proposed \cite{Scheu,Pierer,Pierer4}. Let $ \vm y_i := [\vm s_i\trans, \vm \nu_i\trans, \vm \kappa_i\trans ]\trans \in \mathbb{R}^{p_i}$ be the primal-dual solution of \eqref{eq:local_NLP} and 
	 \begin{equation} \label{eq:local_NLP_Lagrangian}
		\bar{L}_i\indq(\vm y_i) = \bar f_i\indq(\vm s_i) + \vm \nu_i\trans \vm{\bar g}_i\indq(\vm s_i) + \vm \kappa_i\trans \vm{\bar h}_i\indq(\vm s_i)
	 \end{equation}
	the Lagrangian of the local NLPs \eqref{eq:local_NLP} at some iteration~$q$, where $\vm \nu_{i} \in \mathbb R^{n_{g_i}}$ and $\vm \kappa_{i} \in \mathbb R^{n_{h_i}}$ are the local Lagrange multipliers for the constraints \eqref{eq:local_equality} -- \eqref{eq:local_inequality}, respectively. 

The straightforward approach taken in \cite{Scheu}, \cite{Pierer} is to perform the following Newton-like update
\begin{align} \label{eq:NewtonUpdate}
  \vm x_i\indqn = \vm x_i\indq + \vm s_i\indq\,, \quad \vm \lambda_i\indqn = \vm \nu_i\indq\,, \quad \vm \mu_i\indqn = \vm \kappa_i\indq
\end{align}
for each agent $\agents$.
 However, it is shown in \cite{Pierer} that the convergence of this update scheme depends on the coupling strength between subsystems and that the SBDP method might diverge for strongly coupled subsystems since the local steps may be too large to result in a global contraction. A natural extension is to introduce a step size $\alpha>0$ which dampens the updates and prevents large local steps. This modification results in the update law
 \begin{subequations} \label{eq:dampedUpdateLaws}
 \begin{align} \label{eq:dampedNewtonUpdateprimal}
  \vm x_i\indqn = \vm x_i\indq + \alpha \vm s_i\indq 
  \end{align}
  for the primal variable and the update laws
  \begin{align} \label{eq:dampedNewtonUpdatedual}
  \vm \lambda_i\indqn \!=\! \vm \lambda_i \indq \!+\! \alpha(\vm \nu_i\indq - \vm \lambda_i\indq)\,, \smallspace \vm \mu_i\indqn \!=\! \vm \mu_i\indq \!+\! \alpha (\vm \kappa_i\indq - \vm \mu_i\indq )
   \end{align}
  \end{subequations}
for the dual variables. This adaption already guarantees convergence for a broad class of NLPs \eqref{eq:central_NLP}. In particular it can be shown \cite{Pierer4} that the step size $\alpha$ can always be chosen sufficiently small such that the iteration defined by \eqref{eq:dampedUpdateLaws} converges locally if the central NLP is only coupled via the costs and is decoupled in the constraints. However, if we allow for coupled constraints in \eqref{eq:central_NLP}, a more sophisticated update scheme is necessary to ensure convergence. To this end, it is proposed in \cite{Pierer4} that each agent $\agents$ updates its primal and dual variables as follows
\begin{align} \label{eq:primal_dual_update}
\vm p_i\indqn = \vm p_i \indq + \alpha \vm P_i\indq(\vm y_i\indq) (\vm y_i\indq - \vm d_i(\vm p_i\indq))
\end{align}
with offset $\vm d_i(\vm p_i) := [\vm 0\trans, \vm \lambda_i \trans, \vm \mu_i\trans]\trans$ and the matrix-valued function $\vm P_i\indq : \mathbb{R}^{p_i} \rightarrow \mathbb{R}^{p_i \times p_i}$ given as
\begin{align} \label{eq:MixingMatrix}
\vm P_i\indq(\vm y_i) \!=\! \begin{bmatrix}
\nabla_{\vm s_i \vm s_i}^2 \bar{L}_i\indq( \vm y_i) \!&\!\! \nabla_{\vm s_i}\trans \vm{\bar g}_i\indq(\vm s_i) \!\!&\!\! \nabla_{\vm s_i}\trans \vm{\bar h}_i\indq(\vm s_i)\\
-\beta \nabla_{\vm s_i} \vm{\bar g}_i\indq(\vm s_i) \!&\! \vm 0 \!&\! \vm 0 \\
-\beta \vm K_i \nabla_{\vm s_i} \vm{\bar h}_i\indq(\vm s_i) \!\!\!\!\!&\! \vm 0 \!&\! -\beta \vm{\bar H}_i\indq(\vm s_i)
\end{bmatrix}
\end{align}
for each $\agents$. The parameter $\beta >0 $ describes an additional step size for the dual updates, the matrix $\vm K_i:= \diag([\kappa_{1,i},\dots,\kappa_{n_{h_i},i}])$ is a diagonal matrix consisting of each $\kappa_{k,i}$ in $\vm \kappa_i$, $k \in\mathbb{N}_{[1,n_{h_i}]}$, while $\vm{\bar H}_i\indq(\vm s_i)= \diag(\vm{\bar h}_i\indq(\vm s_i))$ is a diagonal matrix consisting of all local inequalities \eqref{eq:local_inequality}. Hereby, $\mathbb{N}_{[1,n_{h_i}]}$ denotes the integer set from $1$ to $n_{h_i}$. The matrix $\vm P_i\indq(\vm y_i)$ defines a structured linear transformation that alters the linearized algorithm to enforce convergence. 

\subsection{Distributed optimization algorithm}
The algorithm in \eqref{eq:NewtonUpdate} is referred to as SBDP \cite{Pierer}, while its transformed update variant \eqref{eq:primal_dual_update} is denoted as SBDP\texttt{+} \cite{Pierer4}. Both are summarized in Algorithm \ref{alg:SBDP} as they only differ in the choice of the update laws \eqref{eq:NewtonUpdate} -- \eqref{eq:primal_dual_update}. The method requires a bi-directional, neighbor-to-neighbor communication network with the same graph structure $\mathcal{G}$ as in the coupling structure of NLP~\eqref{eq:central_NLP}. In Step 1, each agent computes the partial derivative $\nablaxj  L_i\indq$ and shares this quantity with the respective neighbors $\neighs$ for the evaluation of the local cost function \eqref{eq:local_costFunction}. Afterward, each agent solves the NLP \eqref{eq:local_NLP} in parallel to obtain the local primal-dual solution $\vm y_i\indq$. The new central primal-dual solution $\vm p_i\indqn$ is obtained via one of the updates \eqref{eq:NewtonUpdate} -- \eqref{eq:primal_dual_update} in Step~4, before the primal variable $\vm x_i\indqn$ is sent. A stopping criterion is $\|\vm p_i\indqn - \vm p_i\indq\|_\infty \leq \epsilon$ with tolerance $\epsilon>0$. Convergence of SBDP(\texttt{+}) is investigated in \cite{Pierer} and \cite{Pierer4}, respectively, and summarized in the next theorems.
\begin{theorem}[Convergence of SBDP]
Suppose that that LICQ, SOSC, and SCS are satisfied for both problems~\eqref{eq:central_NLP} and~\eqref{eq:local_NLP} and that the generalized diagonal dominance condition, as specified in \cite{Pierer}, holds. 
Then, the iterates generated by the update law~\eqref{eq:NewtonUpdate} converge locally to the KKT-point of NLP~\eqref{eq:central_NLP} at least linearly. 
Moreover, quadratic convergence can be achieved depending on the problem structure.
\end{theorem}
\begin{theorem}[Convergence of SBDP\texttt{+}]
Suppose that the Hessian of the Lagrangian~\eqref{eq:central_Lagrangian} at the optimal solution is positive definite, and that LICQ and SCS are satisfied for both problems~\eqref{eq:central_NLP} and~\eqref{eq:local_NLP}. Then, the iterates generated by the update law~\eqref{eq:primal_dual_update} converge locally to the KKT-point of NLP~\eqref{eq:central_NLP} at a linear rate, provided that the step size~$\alpha$ is sufficiently small and the penalty parameter~$\rho$ is sufficiently large. 
\end{theorem}
\begin{algorithm}[tb] \small
		\caption{ SBDP(\texttt{+}) for solving NLP \eqref{eq:central_NLP} }
		\begin{algorithmic}[1]
    \addtocounter{ALG@line}{-1}
			\State Initialize $\vm p_i^0$; Choose tuning parameters $\alpha,\, \beta,\, \rho$; send  $\vm x_{i}^0$ to all neighbors $\neighs$; set $q\to0$
            \State Compute the mirroring gradient to \eqref{eq:gradient} for all $\neighs$ as
            \begin{equation} \label{eq:gradient_alg}
             \nablaxj  L_i\indq =\nablaxj L_i(\vm x_i\indq, \vm \lambda_i\indq, \vm \mu_i\indq, \Ni{x}\indq)
            \end{equation}
			\State Send $\nablaxj  L_i\indq$ to the respective neighbor $\neighs$.
			\State Solve the local NLP \eqref{eq:local_NLP} to (local) optimality
                \begin{subequations}\label{eq:local_NLP_alg}
                    \begin{alignat}{2}
                        (\vm s_i\indq,\vm \nu_i\indq,\vm \kappa_i \indq)=\arg\min_{\vm s_i} \quad & \bar f_i\indq(\vm s_i) \label{eq:local_costFunctionalg}
                        \\  \st \quad&  \vm  {\bar g}_i\indq(\vm s_i) = \vm 0 \label{eq:local_equalityalg}  \\
                        &\vm {\bar h}_i\indq(\vm s_i)\leq\vm 0\label{eq:local_inequalityalg3}
                    \end{alignat}
                \end{subequations}
            \State For SBDP, update $\vm p_i\indqn$ with \eqref{eq:NewtonUpdate}; For SBDP\texttt{+}, update $\vm p_i\indqn$ either with \eqref{eq:dampedUpdateLaws} or \eqref{eq:primal_dual_update}.
			\State Send $\vm x_i\indqn$ to all neighbors $\neighs$.
			\State Stop if a suitable convergence criterion is met. Otherwise, return to line~$1$ with $q \leftarrow q+1$.
		\end{algorithmic}\label{alg:SBDP}
	\end{algorithm}

\section{Sensitivity-based distributed dynamic optimization}
\label{sec:SBDP_dynamic}
In the context of DMPC, SBDP is used to solve the central OCP in distributed fashion at each sampling step. While it is possible to directly consider discrete-time systems for which the MPC optimization problem takes the form of \eqref{eq:central_NLP}, we will adopt a continuous-time perspective. This viewpoint preserves structural properties in the sense of a ``first optimize then discretize'' approach which is exploited, for instance, by an efficient (gradient-based) solution of the subproblems or for the sensitivity calculation.
\subsection{Problem statement}
In this setting, each agent $\agents$ is governed by the nonlinear, continuous-time system dynamics
	\begin{equation}\label{eq:agent_dynamics}
		\vm{\dot{x}}_i(t) = \vm {f}_i(\vm {x}_i(t), \vm u_i(t), \Ni{x}(t))\,,\smallspace \forall t>0
	\end{equation}
with the states $\vm {x}_i(t)\in \mathbb{R}^{n_{x_i}}$, controls $\vm u_i(t)\in \mathbb{R}^{n_{u_i}}$, initial conditions $\vm {x}_i(0)= \vm x_{i,0}$ and functions $ \vm {f}_i:\mathbb{R}^{n_{x_i}}\times \mathbb{R}^{n_{u_i}} \times \mathbb{R}^{n_{x_{\mathcal{N}_i}}}\rightarrow \mathbb{R}^{n_{x_i}}$. To avoid confusion with optimization variables, state variables are consistently denoted with their explicit time argument, i.e., $\vm x_i(t)$. Similar to \eqref{eq:central_NLP}, the agents solve the central OCP 
\begin{subequations}\label{eq:central_ocp}
 	\begin{alignat}{3}
 		\min_{ \vm {\bar u}(\cdot)} \quad&
 		\sum_{i\in\mathcal V} J_i(\vm {\bar u}_i(\cdot); \Ni{\bar x}(\cdot),\vm x_{i,0})	\label{eq:central_ocp_cost}
 		\\  ~\st \quad&\vm{\dot {\bar x}}_i(\tau) = \vm f_i(\vm {\bar x}_i(\tau), \vm {\bar u}_i(\tau), \Ni{\bar x}(\tau))\,, &&\quad \agents\label{eq:central_ocp_dyn}\\
 		& \vm {\bar x}_i(0) =\vm x_{i,0}\,, &&\quad \agents \label{eq:central_ocp_init}	\\
 		&\vm {\bar u}_i(\tau )\in \mathbb{U}_i\,,\quad \timetau\,, &&\quad \agents\,,\label{eq:central_ocp_constraint}
 	\end{alignat}
 \end{subequations}
 where each input is constrained to the compact set $\mathbb{U}_i\subset \mathbb{R}^{n_{u_i}}$ via the constraint \eqref{eq:central_ocp_constraint} and $\vm{\bar u}(\tau) = [\vm {\bar u}_i(\tau)]_{\agents}$. The bar notation indicates optimization variables depending on the optimization time $\timetau$ in contrast to system variables. 
 The local cost functionals \eqref{eq:central_ocp_cost} are chosen as
\begin{multline} \label{eq:agent_cost}
J_i(\vm {\bar u}_i(\cdot); \Ni{\bar x}(\cdot),\vm x_{i,0}) \\
= V_i(\vm {\bar x}_i(T)) + \int_0^T  l_i(\vm {\bar x}_i(\tau), \vm {\bar u}_i(\tau),\Ni{\bar x}(\tau))\, \dd \tau 
\end{multline}
with prediction horizon $T>0$, integral cost $ l_{i} : \mathbb{R}^{n_{x_i}} \times \mathbb{R}^{n_{u_i}} \times \mathbb{R}^{n_{x_{\mathcal{N}_i}}} \rightarrow \mathbb{R}_{\geq 0}$ and terminal cost  $ V_{i} : \mathbb{R}^{n_{x_i}} \rightarrow \mathbb{R}_{\geq 0}$, $\agents$. The coupling between agents appears via the neighboring states $\vm x_j(t)$ in the dynamics \eqref{eq:agent_dynamics} and the integral cost \eqref{eq:agent_cost}. All functions in OCP~\eqref{eq:central_ocp} are supposed to be sufficiently smooth and we assume the existence and uniqueness of the respective solution $\vm x_i(\cdot; \vm u_i(\cdot),\Ni{x}(\cdot),\vm x_{i,0})$ of \eqref{eq:agent_dynamics} such that the cost functional in~\eqref{eq:agent_cost} is regarded as being dependent on the control trajectory $\vm {\bar u}_i(\cdot)$, external trajectories $\Ni{\bar x}(\cdot)$ and initial condition $\vm x_{i,0}$.

\subsection{Local problems and sensitivities}
Similar to the parametric cased discussed before, modified local OCPs are constructed which are solved at each iteration $q=0,1, \dots$. Specifically, they are given as
\begin{subequations}\label{eq:local_OCP}
	\begin{align}
		\min_{\vm {\bar u}_i(\cdot)} \quad&  \bar J_i\indq(\vm {\bar u}_i(\cdot);\vm x_{i,0}) := J_i(\vm {\bar u}_i(\cdot), \Ni{\bar x}\indq(\cdot),\vm x_{i,0}) \nonumber\\
		&+\sum_{j \in \mathcal{N}_i} \delta J_j(\vm {\bar u}_j\indq(\cdot);\Nj{\bar x}\indq(\cdot),\vm x_{j,0})(\delta \vm{ \bar x}_i(\cdot)) \label{eq:local_OCP_cost}
		\\  ~\st \quad&\vm{\dot {\bar x}}_i(\tau) = \vm f_i(\vm {\bar x}_i(\tau), \vm {\bar u}_i(\tau), \Ni{\bar x}\indq(\tau))\label{eq:local_OCP_dyn}\\
 		& \vm {\bar x}_i(0) =\vm x_{i,0} \label{eq:local_OCP_init}	\\
 		&\vm {\bar u}_i(\tau )\in \mathbb{U}_i\,,\quad \timetau \label{eq:local_OCP_constraint}
	\end{align}
\end{subequations}
The cost $\bar J_i\indq(\vm {\bar u}_i(\cdot);\vm x_{i,0})$ represents the modified local cost functional which is augmented by the sensitivities of the neighboring agents. 
In this continuous-time context, the sensitivity term appearing in the modified local cost function~\eqref{eq:local_OCP_cost}, i.e., $\delta J_j(\vm {\bar u}_j\indq(\cdot);\Nj{\bar x}\indq(\cdot),\vm x_{j,0})(\delta \vm {\bar x}_i(\cdot))$, is defined as the Gâteaux derivative of the neighboring agents' extended cost functional w.r.t.\ the agent states $\vm {\bar x}_i(\tau)$ in direction $\delta \vm {\bar x}_i(\tau) =  \vm {\bar x}_i(\tau) - \vm {\bar x}_i\indq(\tau)$. It is computed as \cite{Pierer3}
\begin{equation}\label{eq:Gateaux_deriv}
\delta J_j(\vm {\bar u}_j(\cdot), \Nj{\bar x}(\cdot))(\delta \vm {\bar x}_i(\cdot))=\!\! \int_0^T \!\! (\vm g_{ji}(\tau))\trans \delta \vm {\bar x}_i(\tau)  \, \dd \tau\,,
\end{equation}
where $\vm g_{ji}(\tau) \in \mathbb{R}^{n_{x_i}}$, $\timetau$ may be interpreted as the time-dependent gradient
\begin{align}
\vm {g}_{ji}(\tau) =&  \nabla_{\vm {x}_i}l_{j}(\vm {\bar x}_j(\tau), \vm {\bar u}_j(\tau),\Nj{\bar x}(\tau)) \nonumber \\ &+ (\nabla_{\vm {\bar x}_i} \vm f_{j}(\vm {\bar x}_j(\tau), \vm {\bar u}_j(\tau),\Nj{\bar x}(\tau))\trans \vm{\bar \lambda}_j(\tau) \label{eq:gradient}
\end{align}
of the cost functional $J_j(\cdot)$, $\neighs$, w.r.t.\ the states $\vm {\bar x}_i(\tau)$, $\agents$. Hereby, $\vm \lambda_j(\tau) \in \mathbb{R}^{n_{x_j}}$ denotes the adjoint state of the neighbors associated with the dynamics  \eqref{eq:local_OCP_dyn}. It may be computed locally in each iteration $q$ via backward integration of the adjoint dynamics
\begin{equation}\label{eq:local_adjdyn}
\vm{\dot {\bar \lambda}}_i(\tau) = - \nablax H_i\indq(\vm {\bar x}_i(\tau), \vm { \bar u}_i(\tau), \vm {\bar \lambda}_i(\tau)) 
\end{equation}
with the terminal condition $\vm {\bar \lambda}_i(T) = \nablax V_i(\vm { \bar x}_i(T))$ and the local Hamiltonian associated with the local OCP \eqref{eq:local_OCP}
\begin{align}\label{eq:local_Hamilton}
H_i\indq(\vm x_i, \vm u_i, \vm \lambda_i) = \,& l_i(\vm x_i, \vm u_i, \Ni{x}\indq) + \vm \lambda_i\trans \vm f_i(\vm x_i, \vm u_i, \Ni{x}\indq) \nonumber \\ &+ \sum_{\neighs}  (\vm g_{ji}\indq)\trans (\vm x_i - \vm x_i\indq)\,,
\end{align}
where the implicit dependency of \eqref{eq:local_Hamilton} on the neighboring states $\vm x_j(\tau)$ and gradients $\vm g_{ji}(\tau)$
is captured by the superscript~$q$. Furthermore, we use the short-hand notation $\vm g_{ji}\indq(\tau)$ to denote the gradient evaluated at iteration $q$. The dynamics \eqref{eq:agent_dynamics} and integral costs in \eqref{eq:agent_cost} are decoupled by treating the neighboring trajectories $\vm {\bar x}_j(\cdot)$ as fixed at the current iteration $q$.
\subsection{Distributed optimal control algorithm}
The decoupled nature of the local OCPs \eqref{eq:local_OCP} is exploited by their parallel solution at the agent level, see Algorithm~\ref{alg:SBDP_cont}. As in Algorithm~\ref{alg:SBDP}, each iteration consists of two computation steps, each followed by neighbor-to-neighbor communication: first exchanging the “mirrored” gradient trajectory~\eqref{eq:gradient}, then broadcasting the current state trajectory. At every iteration~$q$, the local OCP~\eqref{eq:local_OCP} needs to be solved, e.g., via a projected gradient method~\cite{Graichen3} or with the fixed-point iteration in the next section~\cite{Graichen}. Both approaches yield the adjoint state as a byproduct, providing $\vm \lambda_i(\cdot)$ without extra cost. A practical stopping criterion in DMPC is a fixed iteration limit. Convergence of Algorithm~\ref{alg:SBDP_cont} is analyzed in~\cite{Pierer2}.
 \begin{theorem}[Convergence of Algorithm \ref{alg:SBDP_cont}] \label{th:cont_time}
Let $\mathbb{X} \!\subset\! \mathbb{R}^{n_x}$ be a compact set such that $\vm x_0 = [\vm x_{i,0}]_{\agents} \in \mathbb{X}$. 
Assume that OCP~\eqref{eq:central_ocp} admits a unique solution and that both the central and local optimal control laws are locally Lipschitz. 
Then, there exists a sufficiently short prediction horizon $T$ such that Algorithm~\ref{alg:SBDP_cont} converges linearly.
\end{theorem}
Similar to \eqref{eq:dampedUpdateLaws}, the allowable horizon length can be increased by damping the iterates \cite{Pierer2}.
\begin{algorithm}[tb] \small
		\caption{SBDP for solving OCP \eqref{eq:central_ocp}}
		\begin{algorithmic}[1]
      \addtocounter{ALG@line}{-1}
			\State Initialize $\vm {\bar x}_i^0(\tau) = \vm x_{i,0}$, $\vm {\bar \lambda}_i^0(\tau) = \nablax V_i(\vm x_{i,0})$, $\timetau$; set $\vm x_{i,0}$; send $\vm {\bar x}_i^0(\cdot)$ to neighbors $\neighs$; set $q\to0$
            \State Compute the mirroring gradient to \eqref{eq:gradient} for all $\neighs$ as
            \begin{multline}\label{eq:gradientcont_alg}
           \vm {g}_{ij}\indq(\tau) =  \nabla_{\vm { x}_j}l_{i}(\vm {\bar x}_i\indq(\tau), \vm {\bar u}_i\indq(\tau),\Ni{\bar x}\indq(\tau)) \\+ (\nabla_{\vm {x}_j} \vm f_{i}(\vm {\bar x}_i\indq(\tau), \vm {\bar u}_i\indq(\tau),\Ni{\bar x}\indq(\tau))\trans \vm {\bar \lambda}_i\indq (\tau) \,.
            \end{multline}
			\State Send $\vm g_{ij}\indq(\cdot)$ to the respective neighbor $\neighs$.
			\State Compute the trajectories $(\vm {\bar u}_i\indqn(\cdot), \vm {\bar  x}_i\indqn(\cdot), \vm {\bar \lambda}_i\indqn(\cdot))$ by solving the local OCP \eqref{eq:local_OCP}
                \begin{subequations}\label{eq:local_OCP_alg}
                    \begin{alignat}{2}
                        \min_{\vm {\bar u}_i(\cdot)} \quad&  \bar J_i\indq(\vm {\bar u}_i(\cdot);\vm x_{i,0}) \label{eq:local_OCPalg_cost}\\
                      ~\st \quad&\vm{\dot {\bar x}}_i(\tau) = \vm f_i(\vm {\bar x}_i(\tau), \vm {\bar u}_i(\tau), \Ni{\bar x}\indq(\tau))\\
                  & \vm {\bar x}_i(0) =\vm x_{i,0}\label{eq:dist_OCPalg__init}\\
		              &\vm {\bar u}_i(\tau)\in \mathbb{U}_i\,,\quad \timetau \,.\label{eq:local_OCPalg_constraint}
                  \end{alignat}
                \end{subequations}
                \State Send the state trajectory $\vm {\bar x}_i\indqn(\cdot)$ to all neighbors $\neighs$
                \State Stop if a suitable convergence criterion is met. Otherwise, return to line~$1$ with $q \leftarrow q+1$.
		\end{algorithmic}\label{alg:SBDP_cont}
	\end{algorithm}
\subsection{Solution of the local problems via fixed-point iterations}
\label{subsec:fixed_point}
Generally, a suitable local solver is needed to compute the solutions of \eqref{eq:local_OCP} in each iteration $q$. However, for a frequently appearing special case discussed below, the structure of the optimality conditions of \eqref{eq:local_OCP} can be exploited to simplify the solution process. For this purpose, suppose the agent dynamics \eqref{eq:agent_dynamics} are input-affine, i.e., 
\begin{align}\label{eq:agent_dyn_input_affine} 
	\vm f_i(\vm x_i, \vm u_i, \Ni{x}) = \vm f_{i}^0(\vm x_i, \Ni{x}) + \vm B_i(\vm x_i) \vm u_i\,,
\end{align}
with control-independent functions $\vm f_i^0 : \mathbb{R}^{n_{x_i}} \times  \mathbb{R}^{n_{x_{\mathcal{N}_i}}}  \rightarrow \mathbb{R}^{n_{x_i}}$ and matrix functions $\vm B_i : \mathbb{R}^{n_{x_i} }\rightarrow \mathbb{R}^{n_{x_i} \times n_{u_i}}$. Furthermore, the integral costs \eqref{eq:agent_cost} are quadratic in $\vm u_i$, i.e.,
\begin{equation}\label{eq:agent_costs_quadratic}
  l_i(\vm x_i, \vm u_i, \Ni{x}) = l_{i}^0(\vm x_i, \Ni{x}) + \frac{1}{2}\Delta \vm{u}_i\trans \vm R_i  \Delta \vm{u}_i\,,
\end{equation}
with control-independent cost $l_i^0 : \mathbb{R}^{n_{x_i}} \times  \mathbb{R}^{n_{x_{\mathcal{N}_i}}}  \rightarrow \mathbb{R}_{\geq 0}$, diagonal matrix $\vm R_i \succ \vm 0$, and $ \Delta \vm u_i = \vm u_i - \vm u_{i}^{\mathrm{ref}}$ for some reference $\vm u_{i}^{\mathrm{ref}} \in \mathbb{R}^{n_{u_i}}$. 
The constraint set is given as $\mathbb{U}_i =[\vm u_i^-, \vm u_i^+]$ with lower bound $\vm u_i^-$ and upper bound $\vm u_i^+$, respectively. Then, the resulting structure of the optimality conditions of~\eqref{eq:local_OCP} allows for an efficient solution with the fixed-point scheme presented in \cite{Graichen}. Specifically, the first-order optimality conditions for each OCP \eqref{eq:local_OCP} consist of the canonical boundary value
problem (BVP)
\begin{subequations}\label{eq:bvp}
\begin{alignat}{2} \label{eq:sys_dyn}
  \vm{\dot{ \bar x}}_i(\tau) &\!=\! \vm f_i(\vm {\bar x}_i(\tau), \vm {\bar u}_i(\tau), \Ni{\bar  x}\indq(\tau)),
& \vm {\bar  x}_i(0) & \!=\! \vm x_{i,0}
\\ \label{eq:adjoint_dyn}
  \vm {\dot {\bar  \lambda}}_i(\tau) &\!= \!- \nablax H_i\indq(\vm {\bar  x}_i(\tau),\vm {\bar  u}_i(\tau),\vm {\bar  \lambda}_i(\tau)),\,
 & \vm {\bar  \lambda}_i(T) & \!=\! \vm {\bar \lambda}_i^T
\end{alignat}
\end{subequations}
with terminal condition $\vm {\bar \lambda}_i^T = \nablax V_i(\vm {\bar x}_i(T))$, and the pointwise-in-time minimization of the local Hamiltonian~\eqref{eq:local_Hamilton} for all $\timetau$
w.r.t.\  the control, i.e.,\
\begin{equation}\label{eq:PMP}
  \vm {\bar u}_i(\tau) = \argmin_{\vm u_i  \in [\vm u_i^-,\, \vm u_i^+]} 
  H_i\indq(\vm {\bar x}_i(\tau),\vm {u}_i,\vm {\bar \lambda}_i(\tau))\,.
\end{equation}
The problem \eqref{eq:PMP} is strictly convex in $\vm u_i$ 
and  separable for every element $u_{k,i}$ in $\vm u_i$, $k\in \mathbb{N}_{[1,n_{u_i}]}$. Thus, the optimal control $u_{k,i}$ is computed via element-wise projection as
\begin{align} \label{eq:min_control}
  u_{k,i} = \phi_{k,i}(\vm x_i,\vm \lambda_i):= \begin{cases}
     u_{k,i}^-  \smallspace &\text{if } {\tilde u}_{k,i} \leq  u_{k,i}^- \\
     u_{k,i}^+  \smallspace &\text{if } {\tilde u}_{k,i} \geq  u_{k,i}^+ \\
    {\tilde u}_{k,i}  \smallspace &\text{if } {\tilde u}_{k,i} \in ( u_{k,i}^-, u_{k,i}^+)
  \end{cases}
\end{align}
with the (vector-valued) unconstrained minimizer 
\begin{equation} \label{eq:unconstrained_minimizer}
	\vm {\tilde u}_{i} = \vm u_{i}^{\mathrm{ref}}- \vm R_i^{-1} \vm B_i\trans(\vm x_i) \vm \lambda_i\,,
\end{equation}
where the element-wise control functions \eqref{eq:min_control} are summarized as $\vm \phi_i(\vm x_i,\vm \lambda_i) := [\phi_{k,i}(\vm x_i,\vm \lambda_i)]_{k \in \mathbb{N}_{[1,n_{u_i}]}}$. 
By inserting \eqref{eq:min_control} into the canonical equations
\eqref{eq:bvp}, we obtain the following BVP 
\begin{subequations}\label{eq:can}
\begin{alignat}{2}
    \vm{\dot{ \bar  x}}_i(\tau) &= \vm F_i\indq(\vm {\bar  x}_i(\tau),\vm { \bar  \lambda}_i(\tau))\,, \quad &\vm { \bar  x}_i(0) &= \vm x_{i,0} \label{eq:can_1} \\
    \vm{\dot{ \bar \lambda}}_i(\tau) &= \vm G_i\indq(\vm {\bar  x}_i(\tau), \vm { \bar  \lambda}_i(\tau))\,, \quad &\vm { \bar  \lambda}_i(T) &= \vm {\bar \lambda}_{i}^T \label{eq:can_2}
  \end{alignat}
\end{subequations}
with the functions  $\vm F_i\indq(\vm x_i, \vm
\lambda_i):= \vm f_i(\vm x_i, \vm \phi_i(\vm x_i, \vm
\lambda_i), \Ni{x}\indq)$ and  $\vm G_i\indq(\vm x_i, \vm \lambda_i):= - \nablax
H_i\indq (\vm x_i,  \vm \phi_i(\vm x_i, \vm \lambda_i), \vm \lambda_i)$, where the superscript captures the implicit dependency on external trajectories at iteration $q$ of both functions. 

The separation of the boundary conditions is exploited in the numerical solution of the BVP \eqref{eq:can} with the fixed-point scheme in Algorithm \ref{alg:fixed_point}.
\begin{algorithm}[tb] \small 
	\caption{Fixed-point iterations for solving \eqref{eq:local_OCP}}
	\begin{algorithmic}[1]
    \addtocounter{ALG@line}{-1}
		\Statex Initialize $ \vm{ \bar \lambda}_i^{0|q}(\tau) = \vm {\bar \lambda}_i^{q}(\tau) $; choose $j_{\max}$; set $j \to 1$ 
		\State \textbf{While} $ j \leq j_{\max}$ \textbf{do}
		\State \hskip 1em Compute $\vm {\bar x}_i\indjq(\tau)$ via forward integration of
		\begin{align} \label{eq:int_x}
			\vm{\dot {\bar  x}}_i\indjq(\tau)&= \vm F_i\indq(\vm {\bar  x}_i\indjq(\tau), \vm {\bar \lambda}_i\indjpq(\tau))\,, \quad \vm {\bar  x}_i\indjq(0) = \vm x_{i,0}\,.
		\end{align} 
		\State \hskip 1em Compute $\vm \lambda_i\indjq(\tau)$ via backward integration of 
		\begin{align} \label{eq:int_lambda}
			\vm{\dot {\bar \lambda}}_i\indjq(\tau) = \vm G_i\indq(\vm {\bar  x}_i\indjq(\tau), \vm {\bar \lambda}_i\indjq(\tau))\,, \quad  \vm {\bar \lambda}_i^{j|k}(T) = \vm{\bar \lambda}_i^T
		\end{align}
		\hskip 1em with terminal condition $\vm{\bar \lambda}_i^T = \nablax V_i(\vm {\bar x}_i\indjq(T))$.
		\State \hskip 1em Set $j \leftarrow j+1$.
		\State \textbf{End While}
		\State Compute $\vm {\bar u}_i\indqn(\tau)=\vm \phi_i(\vm {\bar x}_i^{j_{\max}|q}(\tau),\vm {\bar \lambda}^{j_{\max}|q}_i(\tau))$, set $\vm {\bar x}_i\indqn(\tau) =\vm {\bar x}_i^{j_{\max}|q}(\tau)$, $\vm {\bar \lambda}_i\indqn(\tau)=\vm {\bar\lambda}_i^{j_{\max}|q}(\tau)$, $\timetau$, and return to Step 4 of Algorithm \ref{alg:SBDP_cont}.
	\end{algorithmic}\label{alg:fixed_point}
\end{algorithm}
The scheme is easy to implement, since
the fixed-point iteration merely consists of the sequential
forward and backward integration of \eqref{eq:can_1} -- \eqref{eq:can_2}. The algorithm requires an initial guess for the local adjoint state trajectory $\vm \lambda_i(\cdot)$ which is set in a warm start fashion to the value of the corresponding previous iterate from Algorithm \ref{alg:SBDP_cont}. 
In a practical implementation, a discretization
of the time grid $[0,\,T]$ is required for the numerical integration.
With regard to a real-time capable DMPC scheme, integration methods with a fixed step size are suitable, whereby the number of discretization points affects the optimality, computation time and communication effort. Furthermore, the adjoint state trajectory $\vm \lambda_i(\cdot)$ is part of the solution and can be readily used in the evaluation of the gradient~\eqref{eq:gradient}. The convergence of Algorithm \ref{alg:SBDP_cont} under inexact minimization of the local OCPs \eqref{eq:local_OCP} with Algorithm \ref{alg:fixed_point} is investigated in \cite{Pierer3}.

\begin{theorem}[Convergence under Inexact Minimization]
Consider Algorithm~\ref{alg:SBDP_cont} applied to a neighbor-affine form \cite{Burk} of the problem class specified by \eqref{eq:agent_dyn_input_affine} -- \eqref{eq:agent_costs_quadratic}, where the OCPs~\eqref{eq:local_OCP} are solved inexactly with a fixed number of iterations $j_{\max}$ of Algorithm \ref{alg:fixed_point}. Under the assumptions of Theorem~\ref{th:cont_time}, there exists a sufficiently short prediction horizon $T(j_{\max})$ ensuring linear convergence of Algorithm~\ref{alg:SBDP_cont}.
\end{theorem}

\section{Application to distributed model predictive control} 
\label{sec:DMPC}
In this section, it is discussed how the distributed solution of the central OCP \eqref{eq:central_ocp} via Algorithm \ref{alg:SBDP_cont} is employed in a real-time capable DMPC scheme. To this end, we first briefly recap a central MPC scheme without terminal constraints on which the distributed approach is based. Then, the sensitivity-based DMPC scheme is presented.  
\subsection{Central MPC scheme}
A centralized MPC scheme for setpoint stabilization relies on the repeated online solution of the central OCP \eqref{eq:central_ocp} at each sampling point $t_k = k \Delta t$, $k \in \mathbb{N}_{0}$ with the currently measured (or estimated) system state $\vm {\bar x}_i(0) = \vm x_i^k = \vm x_i(t_k)$ to control each subsystem \eqref{eq:agent_dynamics} to a predefined reference point $(\vm x_{\mathrm{ref}},\vm u_{\mathrm{ref}})$ with $\vm f(\vm x_{\mathrm{ref}},\vm u_{\mathrm{ref}})=\vm 0$, $\vm f(\cdot) = [\vm f_i(\cdot)]_{\agents}$. Usually, MPC strategies assume that the optimal solution of the OCP \eqref{eq:central_ocp} is known at each $t_k$, such that the first part of the optimal control $\vm{\bar u}_i^{*|k}(\tau)$, $\timetau$, i.e., 
\begin{align} \label{eq:optimal_control_law}
  \vm u_i(t_k + \tau ) = \vm {\bar u}_i^{*|k}(\tau)\,, \quad \timetausample\,, \quad \agents 
\end{align}
is applied to each system \eqref{eq:agent_dynamics}. In the next sampling step $t_{k+1} = t_k + \Delta t$, the problem \eqref{eq:central_ocp} is solved again with the new initial value $\vm x_i^{k+1}$. In the nominal case, the next system state is given by the predicted value $\vm x_i^{k+1} = \vm {\bar x}_i^{*|k}(\Delta t)$, where $\vm {\bar x}_i^{*|k}(\cdot; \vm {\bar u}_i^{*|k}(\cdot), \Ni{\bar x}^{*|k}(\cdot), \vm x_{i,0})$ is the optimal state trajectory in MPC step $k$. The OCP \eqref{eq:central_ocp} is formulated on purpose without terminal constraints to reduce the computational burden and to possibly allow for the efficient solution of the local problems with Algorithm \ref{alg:fixed_point}. In order to formally guarantee the stability of an MPC scheme without terminal constraints, additional requirements on the central OCP \eqref{eq:central_ocp} are necessary. In most cases, it is assumed that there exists a local feedback law $\vm u = [\vm r_i(\vm x)]_{\agents}\in \mathbb{U}$, $\mathbb{U} := \mathbb{U}_1\times \dots \times \mathbb{U}_M$ and that the terminal cost $V (\vm x)= \sum_{\agents} V_i(\vm x_i)$ constitutes a local control Lyapunov function (CLF) on a compact, invariant set $\Omega_\beta = \{\vm x \in \mathbb{R}^n \,|\, V(\vm x) \leq \beta\} $ \cite{Mayne,Chen,Limon,Graichen2}. Then, asymptotic stability of the closed loop under the MPC control law \eqref{eq:optimal_control_law} can be shown for a certain domain of attraction $\Gamma \subset \mathbb{R}^n$, i.e., for all $ \vm x_0 \in \Gamma$, the MPC controller asymptotically stabilizes the reference point $\vm x_{\mathrm{ref}}$ of the nominal system~\eqref{eq:agent_dynamics} \cite{Limon, Graichen, Graichen2}.
\subsection{Real-time DMPC approach}
In order to use SBDP in the context of real-time DMPC, only a fixed number of iterations $q_{\mathrm{max}}$ of Algorithm \ref{alg:SBDP_cont} can be performed per sampling step. Terminating the algorithm after $q_{\max}$ iterations yields the suboptimal control $\vm {\bar u}_i^{q_{\max}| k}(\tau)$, $\timetau$ of which the first part is applied to the system \eqref{eq:agent_dynamics}, i.e., 
\begin{align} \label{eq:suboptimal_control_law}
  \vm u_i(t_k + \tau ) = \vm {\bar u}_i^{q_{\max}| k}(\tau)\,, \quad \timetausample\,.
\end{align}
In the next DMPC step, Algorithm \ref{alg:SBDP_cont} is re-initialized with the state and adjoint state trajectories of the previous DMPC step as a warm start. These trajectories can be shifted by $\Delta t$, where the undefined part for $\tau \in (T,\, T + \Delta t]$ may be obtained with the terminal control law $\vm r(\vm x(\tau))$. The DMPC scheme is summarized in Algorithm~\ref{alg:real_time_DMPC}. 

\begin{algorithm}[tb]\small
	\caption{Real-time sensitivity-based DMPC}
  \begin{algorithmic}[1] 
  \addtocounter{ALG@line}{-1}
	\State Initialize $\vm{\bar x}_i^{0|0}(\tau)$, $\vm{\bar \lambda}_i^{0|0}(\tau)$, $\timetau$; choose $q_{\max}$; set $k\to 0 $
		\State Obtain state measurement $\vm x^k_{i} = \vm x_i(t_k)$.
		\State Initialize Algorithm \ref{alg:SBDP_cont} with $\vm {\bar x}_i^{0|k}(0) = \vm x_i^k$,  $\vm {\bar x}_i^{0|k}(\cdot)$ and $\vm {\bar \lambda}_i^{0|k}(\cdot)$, and perform $q_{\mathrm{max}}$ iterations to obtain the trajectories $(\vm {\bar u}_i^{q_{\max}| k}(\tau), \vm {\bar x}_i^{q_{\max}| k}(\tau), \vm {\bar\lambda}_i^{q_{\max}| k}(\tau))$, $\timetau$.
		\State Apply $\vm {\bar u}_i^{q_{\max}| k}(\tau)$, $\timetausample$ to the system \eqref{eq:agent_dynamics}.
		\State Update $\vm {\bar x}_i^{0|k+1}(\tau) = \vm {\bar x}_i^{q_{\max} |k}(\tau) $, $\vm {\bar \lambda}_i^{0| k+1}(\tau) = \vm {\bar\lambda}_i^{q_{\max}|k}(\tau) $, $\timetau$, and go to Step 1 with $k\leftarrow k+1$.
	\end{algorithmic} \label{alg:real_time_DMPC}
\end{algorithm}

In the first DMPC step $k=0$, an initial guess for the state and adjoint state trajectories $\vm {\bar x}_i^{0|0}(\tau)$ and $\vm {\bar \lambda}_i^{0|0}(\tau)$, $\timetau$ is needed for the initialization of Algorithm~\ref{alg:SBDP_cont}. A possible way to determine these trajectories in a distributed fashion is by forward integration of
\begin{equation}
	 \vm{\dot {\bar x}}_i^{0|0}(\tau) = \vm f_i(\vm {\bar x}_i^{0|0}(\tau), \vm {u }_{i}^{\mathrm{ref}}, \vm {x}_{\mathcal{N}_i}^{ \mathrm{ref}})
\end{equation}
with $\vm x_i^{0|0}(0) = \vm x_{i,0}$ and local references $[\vm x_i^{\mathrm{ref}}]_{\agents} = \vm x_{\mathrm{ref}}$, $[\vm u_i^{\mathrm{ref}}]_{\agents} = \vm u_{\mathrm{ref}}$, followed by a backward integration of
\begin{equation}\label{eq:}
\vm{\dot {\bar\lambda}}_i^{0|0}(\tau) = - \nablax H_i^{0|0}(\vm {\bar x}_i^{0|0}(\tau), \vm u_{i}^{\mathrm{ref}}, \vm {\bar \lambda}_i^{0|0}(\tau)) 
\end{equation}
with the terminal condition $\vm {\bar \lambda}_i^{0|0}(T) = \nablax V_i(\vm{ \bar x}_i^{0|0}(T))$. Furthermore, we set $\Ni{\bar x}^{0|0}(\tau) = [\vm x_{j}^{\mathrm{ref}}]_{\neighs}$ and $\vm g_{ji}^{0|0}(\tau) = \vm 0$, $\timetau$, $\neighs$ to resolve the implicit dependencies of the Hamiltonian \eqref{eq:local_Hamilton}. Hereby, only the local state reference values $\vm x_i^{\mathrm{ref}}$ need to be sent offline to the neighbors. 
In \cite{Pierer3}, the stability properties of Algorithm \ref{alg:real_time_DMPC} are investigated.
\begin{theorem}[Real-time DMPC stability]
Under the assumptions that Algorithm \ref{alg:SBDP_cont} is convergent by a suitable choice of the prediction horizon $T$, the number of iterations $q_{\mathrm{max}}$ is sufficiently large, and that the initial optimization error is sufficiently small, the reference of the closed loop system under the control law \eqref{eq:suboptimal_control_law} is exponentially stable and the optimization error decreases incrementally.
\end{theorem}
\section{DMPC of a water tank system}
\label{sec:Experiment}
To demonstrate the practical relevance of the real-time sensitivity-based approach, the setpoint stabilization of a coupled water tank system is considered. We first describe the experimental setup and utilize Algorithm \ref{alg:SBDP} to estimate unknown model parameters in distributed fashion. Then, the sensitivity-based DMPC controller is designed to guarantee nominal stability and evaluated in numeric simulations. Finally, experimental results of the sensitivity-based DMPC controller for a setpoint change and different disturbance scenarios are presented. 
\subsection{Setup and parameter estimation}
\begin{figure}[tb]
		\centering
		\includegraphics{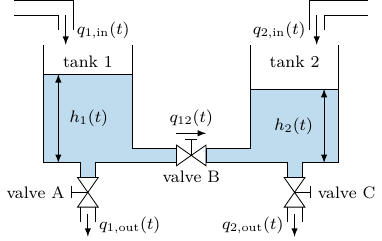}
		\caption{Schematic diagram of the the coupled tank system.}
		\label{fig:watertank_schematic}
		\hspace{1mm}
\end{figure}

\begin{figure}[tb]
		\centering
		\includegraphics{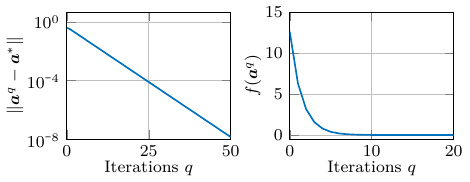}
		\caption{Convergence the error (left) and value function (right) of Algorithm \ref{alg:SBDP} applied to the parameter estimation problem \eqref{eq:parameter_QP}.}
		\label{fig:Parameter_estimation}
\end{figure}

The coupled water tank setup consists of two interconnected tanks, each subject to external inflows or outflows, see Figure \ref{fig:watertank_schematic}. The nonlinear dynamics of each tank $ i \in \mathcal{V}=\{1,\,2\}$, $\mathcal{N}_i = \mathcal{V} \setminus \{i\}$ follow from the mass balance of inflows and outflows for a constant water density as 
\begin{equation} \label{eq:model_watertanks}
	\dot{h}_i(t) = q_{i,\mathrm{in}}(t) + q_{i,\mathrm{out}}(h_i(t)) + \sum_{\neighs} q_{ij}(h_i(t), h_j(t))\,,
\end{equation}
where the water height $h_i(t) = x_i(t) \geq 0 $ represents the state of each agent and is measured by differential pressure sensors. 
The individual functions $ q_{i,\mathrm{in}} = \frac{1}{A}u_i $, $q_{i,\mathrm{out}}(h_i) = \frac{a_{i}}{A}\sqrt{2gh_i}$ and $q_{ij}(h_i, h_j) = \frac{a_{ij}}{A} \sign(h_j - h_i)\sqrt{2g|h_i - h_j|}$ describe the inflow and outflow of each tank, respectively, as well as the coupling flow between adjacent tanks. Here, $A$ is the base area of the tanks, $a_i$, $a_{ij}$ are the effective cross section areas of the outflow and coupling flow respectively, and $g$ is the gravitational constant. Since the coupling function is non-differentiable at $h_i = h_j$, its derivative is approximated by a third-order polynomial for $|h_i - h_j| \leq 0.5\, \si{\centi\meter}$. The volume flows $u_i(t)$, provided by two pumps, act as an input for each tank and are constrained to the interval $[u_i^-,\,u_i^+]$. The nominal configuration of the coupled tank system which is also used as the DMPC internal prediction model is: valve A closed, valve B open and valve C open. 
\tikzexternalenable
\tikzexternaldisable
The values of the effective cross section areas are determined experimentally by steering the nominal system to different equilibrium points and computing $a_2$, $a_{12}$, and $a_{21}$ via a least-square approach
\begin{subequations}\label{eq:parameter_QP}
	\begin{alignat}{2}
		\min_{\vm a} &\quad  f(\vm a) = \|\vm M \vm a - \vm y \|^2 \label{eq:parameter_costFunction}\\
		~\st &\quad  a_{12} = a_{21} \label{eq:constraint_parameter_QP}\\
		 &\quad a_2,\, a_{12},\, a_{21}\geq 0 \label{eq:inequality_constraint_parameter_QP}
	\end{alignat}
\end{subequations}
with regressor $\vm a = [a_{12},a_2,a_{21}]\trans$, coefficient matrix $\vm M \in \mathbb{R}^{m \times n}$ and observations $\vm y \in \mathbb{R}^{m \times 1}$. Overall $m=100$ data points are recorded. The constraints \eqref{eq:constraint_parameter_QP} -- \eqref{eq:inequality_constraint_parameter_QP} arise from physical considerations. If we assign $a_{12}$ with its positivity constraints to agent $1$ and $a_2$, $a_{21}$ with the their positivity constraints together with the equality constraint \eqref{eq:constraint_parameter_QP} to agent $2$, then problem \eqref{eq:parameter_QP} takes the form of~\eqref{eq:central_NLP}. Hence, Algorithm \ref{alg:SBDP} may be applied to estimate the parameter values from the measurement data. Hereby, we additionally substitute the equality constraint \eqref{eq:constraint_parameter_QP} into the cost function of agent 2 to simplify its problem. Figure \ref{fig:Parameter_estimation} shows the convergence of Algorithm \ref{alg:SBDP} with the update rule \eqref{eq:NewtonUpdate} for an initial guess of $a_1=a_{12}=a_{21}=0$ toward the centrally obtained solution~$\vm a\inds$. After $20$ iterations the optimal value function evaluates to $f(\vm a^{20}) = 1.575\times 10^{-4} $  showing a successful distributed estimation of the parameters.
\tikzexternalenable
\tikzexternaldisable
The resulting values of the effective cross sections together with other relevant parameters are summarized in Table~\ref{tab:parameters}. 
\begin{table}[tb] 
\centering
\begin{tabularx}{\columnwidth}{>{\centering\arraybackslash}l|*{6}{|>{\centering\arraybackslash}X}}
Symbol & $A$ & $a_i$&  $a_{ij}$ & $g$ & $ u_i^{-}$ & $u_i^{+}$\\ 
Value & $144$ & $0.354$ & $0.216$ & $981$ & $8.333$ & $100$\\ 
Unit & $\si{\centi\meter\squared}$ & $\si{\centi\meter\squared}$ & $\si{\centi\meter\squared}$ & $\si{\centi\meter\per\second\squared}$& $\si{\centi\meter\cubed \per \second}$&$\si{\centi\meter\cubed \per \second}$\\ 
\end{tabularx}
\caption{Parameter values of the coupled tank system.}\label{tab:parameters}
\end{table}
\subsection{DMPC design and evaluation}
For the control task, we consider the stabilization of the reference point $\vm x_{\mathrm{ref}} = [40 \,\si{\centi \meter},\, 20\,\si{\centi \meter}]\trans$ and $\vm u_{\mathrm{ref}} = [44.27\, \si{\centi\meter\cubed\per \second}\,, 27.24\, \si{\centi\meter\cubed\per \second}] $.
To this end, we choose a quadratic integral cost function  
\begin{align} \label{eq:costs_watertanks}
	l_i(\vm x_i, \vm u_i, \Ni{x}) 
	&= Q_i \Delta x_i^2 + R_i \Delta u_i^2
\end{align}
with $ \Delta x_i := x_i - x_{i}^{\mathrm{ref}}$, $ \Delta u_i := u_i - u_{i}^{\mathrm{ref}}$, $Q_i= 1\,\si{\per\centi\meter\squared}$, and $R_i= 0.1\,\si{\second\squared\per\centi\meter^{6}}$. We design a separable terminal cost $V_i(x_i) = P_i \Delta x_i^2$, and terminal control law $r_i(\vm x) = -K_i \Delta x_i - \sum_{\neighs} K_{ij} \Delta x_j + u_i^{\mathrm{ref}}$, via the LMI-based approach in \cite{Pierer3} which results in $P_1 =48.30\, \si{\per\centi\meter\squared} $, $P_2 =30.87\, \si{\per\centi\meter\squared}$, $K_1 = 3.06$, $K_2 = 1.97$, $K_{12} = K_{21} = 0$, where $\gamma = 1.075$ in the notation of \cite{Pierer3} is used. The terminal region in which the constraints and the CLF inequality are satisfied is computed numerically via a sampling-based approach and follows as $\Omega_\beta =\{\vm x \in \mathbb{R}^2\, |\, V(\vm x) \leq 6.334\times 10^{3}\}$ \cite{Chen}. We set the prediction horizon to $T=6\,\si{\second}$ and evaluate the stabilization of the reference point from a given initial state $\vm x_0$. Rather than checking if the initial state is within the theoretical region of attraction $\Gamma$ given by the (conservative) estimate in \cite{Limon,Graichen2}, we numerically compute $\Gamma$ by checking if the predicted central terminal state $\vm{\bar x}^{*|k}(T) = [\vm x_i^{*|k}(T)]_{\agents}$ reaches the terminal region in each MPC step $k$, i.e., $\vm{\bar x}^{*|k}(T) \in \Omega_\beta$. For the stabilization of $\vm x_{\mathrm{ref}}$ from the initial state $\vm x_0 = [30\,\si{\centi\meter}, 35 \,\si{\centi\meter}]\trans$, we have $V(\vm{\bar x}^{*|k}(T))\leq 5.869 \times 10^3 \leq \beta $ for all $k \in \mathbb{N}_0$ which implies $\vm{\bar x}^{*|k}(T) \in \Omega_\beta $  and thus $\vm x_0 \in \Gamma$.  Therefore, the central, optimal MPC scheme is asymptotically stable for this particular initial state and the real-time DMPC stability considerations in \cite{Pierer2} apply. 

For the evaluation of the real-time sensitivity-based controller, we implement Algorithm \ref{alg:SBDP_cont}, where the local OCPs \eqref{eq:local_OCP_alg} are solved with the fixed-point iteration scheme of Algorithm \ref{alg:fixed_point}, in C\texttt{++} code with a \textsc{Matlab} Cmex-interface. Algorithm \ref{alg:fixed_point} is applicable since the dynamics \eqref{eq:model_watertanks}, costs \eqref{eq:costs_watertanks}, and constraints are in the form as required in Section \ref{subsec:fixed_point}. The sampling time is set to  $\Delta t = 200\, \si{\milli\second}$ and the prediction interval $[0,T]$ is discretized into $30$ equidistant subintervals. A Heun integration scheme is used for the numerical integration of \eqref{eq:int_x} -- \eqref{eq:int_lambda}.  

At first, we investigate how different combinations of inner ($j_{\max}$) and outer ($q_{\max}$) iterations affect the control quality of the sensitivity-based DMPC controller. 
As measure of the optimality, we inspect the closed loop cost 
\begin{equation}
	J_{cl} = \frac{1}{T_{\mathrm{sim}}}\sum_{\agents}\int_{0}^{T_{\mathrm{sim}}} Q_i \Delta x_i^2(t) + R_i  \Delta u_i^2(t)\, \dd t
\end{equation}
which is calculated with the closed loop trajectories $x_i(\cdot)$ and $u_i(\cdot)$ resulting from applying the DMPC control law for a simulation time of $T_{\mathrm{sim}} = 150\, \si{\second}$. Table \ref{tab:optimality} shows the value of $J_{cl} $ for different combinations of $q_{\max}$ and $j_{\max}$. For $q_{\max}$, $j_{\max}\geq 3$, the difference between distributed and central control is less than $3 \times 10^{-2}$, which illustrates that already a very low number of inner and outer iterations lead to a practically optimal MPC behavior. We set $q_{\max} =3$ and $ j_{\max} = 5$ for the following considerations.
\begin{table}[tb] 
\centering
\begin{tabularx}{\columnwidth}{>{\centering\arraybackslash\normalfont}c|*{3}{>{\centering\arraybackslash\normalfont}X}}
\begin{tikzpicture}[scale = 0.5,baseline=(current bounding box.center)]
  \draw (0,1) -- (0.8,0.2); 
  \node[anchor=south east] at (0,0) {$q_{\max}$}; 
  \node[anchor=north west] at (1,1) {$j_{\max}$}; 
\end{tikzpicture}& 1 & 3 & 5 \\ \hline
1 & 64.487 & 53.366& 51.700\\ \hline
3 & 54.061 & 51.307 & 51.282  \\ \hline
5 & 51.876  & 51.282   & 51.282   \\ 
\end{tabularx}
\caption{Closed loop cost $J_{cl}$ for different combinations of inner ($j_{\max}$) and outer iterations ($q_{\max}$) of the DMPC approach. The value obtained from the MPC simulation is $J_{cl}= 51.282$.}\label{tab:optimality}
\end{table}
Since the only computations performed locally are the evaluation of the sensitivities in form of the gradients \eqref{eq:gradient} and the forward/backward integration of the canonical equations~\eqref{eq:can}, the computational footprint at the agent level is minimal. In fact, for $q_{\max} =3$ and $ j_{\max} = 5$, Algorithm \ref{alg:SBDP_cont} takes around $250\, \si{\micro\second}$ per agent on an Intel Core i7 ($2.80\, \si{\giga \hertz}$), excluding the communication overhead. Regarding the communication effort, $180$ floats must be exchanged between the controllers at every DMPC sampling step. Depending on the communication protocol and infrastructure, this data transmission may take up a substantial amount of the overall execution time available in each sampling step \cite{Burk,Stomberg2,Pierer5}. However, for this simulation the communication time is negligible since Algorithm 2 is executed in a non-distributed environment.

\begin{figure}[tb]
		\centering
		\includegraphics{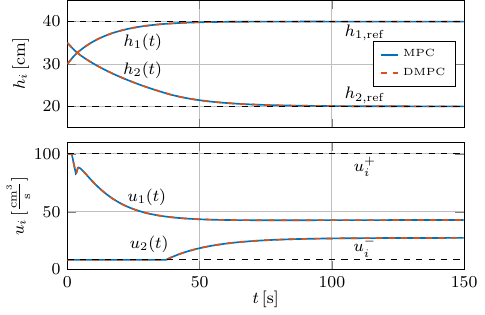}
		\vspace{-3mm}
		\caption{Simulative 
		comparison of the closed loop (D)MPC trajectories of the multi-tank setup for $q_{\max} = 3$ and $ j_{\max} = 5$.}
		\label{fig:Simulation_Watertank}
\end{figure}

Figure \ref{fig:Simulation_Watertank} shows the resulting closed loop trajectories resulting from an optimal MPC scheme where the central OCP \eqref{eq:central_ocp} is solved to optimality at each sampling point $t_k$ and the sensitivity-based DMPC scheme described in Algorithm \ref{alg:real_time_DMPC}. The MPC and DMPC control inputs exhibit nearly no discrepancies, even at the crossing point of the water levels $h_1$ and $h_2$ at approximately $5\, \si{\second}$, where the nonlinearity and the approximation of the derivative in the model \eqref{eq:model_watertanks} becomes significant. This underscores that the employed sensitivity-based DMPC approach is able to achieve practically the same control performance as the centralized scheme with a very limited computational effort at the agent level. 
\tikzexternalenable
\tikzexternaldisable

\subsection{Experimental results}
\tikzexternalenable
\begin{figure}[tb]
		\centering
		\includegraphics{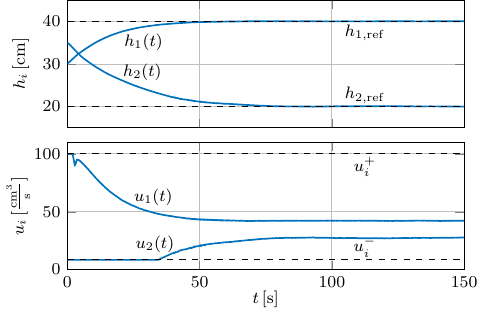}
		\vspace{-4mm}
		\caption{Experimental results of the sensitivity-based DMPC scheme applied to the coupled water tank setup.}
		\label{fig:Experiment_Watertank}
\end{figure}
\tikzexternaldisable
In addition to the simulation study presented in the previous section, the applicability of the sensitivity-based DMPC scheme is demonstrated via experiments on the real multi-tank setup. To this end, Algorithm~\ref{alg:real_time_DMPC}, in which the local OCPs are solved with Algorithm~\ref{alg:fixed_point}, is implemented on the corresponding d\textsc{SPACE} hardware via \textsc{Matlab/Simulink}. Due to this hardware restriction, Algorithm~\ref{alg:real_time_DMPC} is only executed in a pseudo-distributed fashion in the sense that the data separability and an independent execution of the computation steps are ensured, but no actual communication is realized. Hereby, all parameters are set as in the simulation.
\begin{table}[tb] 
\centering
\begin{tabularx}{\columnwidth}{>{\centering\arraybackslash\normalfont}l|*{3}{>{\centering\arraybackslash\normalfont}X}}
& Valve A & Valve B & Valve C \\ \hline
Nominal & closed & open & open \\ \hline
Scenario $\mathsf{I}.$ & open & open& open\\ \hline
Scenario $\mathsf{II}.$ & closed & open & closed  \\ \hline
Scenario $\mathsf{III}.$ & closed  & closed   & open   \\ 
\end{tabularx}
\caption{Nominal configuration and different disturbance scenarios.}\label{tab:disturbances}
\end{table}
\tikzexternalenable
\begin{figure}[tb]
		\centering
		\vspace{5mm}
		\includegraphics{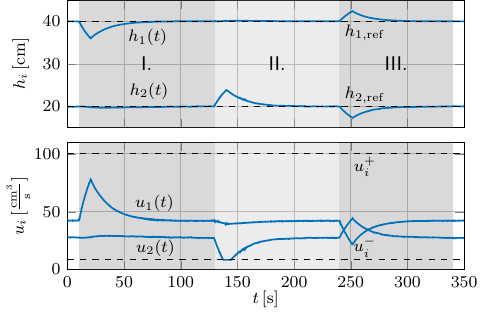}
		\vspace{-4mm}
		\caption{Experimental results of the sensitivity-based DMPC controller for different disturbance scenarios specified in Table \ref{tab:disturbances}.}
		\label{fig:Experiment_Watertank_disturbance}
\end{figure}
\tikzexternaldisable
The experimental results of the DMPC controller, shown in Figure \ref{fig:Experiment_Watertank}, exhibit good agreement with the simulation results in Figure \ref{fig:Simulation_Watertank}. This is underpinned by the closed loop cost value of $J_{\mathrm{cl}} = 54.633$ which is about $6.5\, \si{\percent}$ higher than the value from the simulations. This difference is primarily attributed to nonlinear pump effects which are unaccounted for.
From a practical viewpoint it is also of interest how the sensitivity-based DMPC controller rejects disturbances. To this end, we subject the multi-tank setup to three different disturbance scenarios $\mathsf{I}.$ to $\mathsf{III}.$ which together with the nominal configuration are summarized in Table~\ref{tab:disturbances}. Each scenario regards a different configuration of the three valves A-C, see Figure \ref{fig:watertank_schematic}. The corresponding experimental results are depicted in Figure \ref{fig:Experiment_Watertank_disturbance}. Hereby, the system is initially at the reference point $(\vm x_{\mathrm{ref}},\, \vm u_{\mathrm{ref}})$ before the different disturbance scenarios $\mathsf{I}.$ to $\mathsf{III}.$ are applied to the water tanks, where each valve is manually opened for a duration of $10\, \si{\second}$. After the system has settled at the reference, the next scenario is  started. The duration of each of the scenarios $\mathsf{I}.$ to $\mathsf{III}.$ is indicated by the shaded areas in Figure \ref{fig:Experiment_Watertank_disturbance}. Although the DMPC internal model does not account for these disturbances, the distributed controller is able to successfully reject them. 

\section{Conclusion}
\label{sec:Conclusion}
This paper discusses a sensitivity-based framework for large-scale distributed optimization. By augmenting local subproblems with first-order sensitivity information from neighboring agents and a primal decomposition approach, the method provides formal convergence guarantees to the optimal solution while requiring only local computation and neighbor-to-neighbor communication. Different update rules are introduced such as Newton-like steps for fast progress, damped variants that have a stabilizing effect, and the SBDP\texttt{+} update scheme that provides local convergence for all coupling structures. 

The approach naturally extends to dynamic optimization and distributed MPC, where local gradients can be obtained efficiently via adjoint-based integration. Hereby, the convergence is linked to the prediction horizon length. The efficacy of the method is demonstrated with the real-time DMPC of coupled water tanks with guaranteed stability. The DMPC scheme achieved a similar control quality as the central MPC with only a few SBDP iterations. 

Regarding methodology, future work will address global convergence guarantees, application to NLPs with non-differentiable objectives and more detailed convergence analyses in continuous-time. On the practical side, learning-based approximation of the gradients, adaptive step sizes and application to large scale machine learning problems are of interest.  

\printbibliography

\end{document}

%% file: commands.tex
\usepackage{amsmath,amsfonts,amssymb}
\usepackage{dsfont}
\usepackage{algorithm}
\usepackage{algpseudocode}
\usepackage[per-mode=fraction]{siunitx}
\usepackage{tikz}
\usepackage{pgfplots}
\usepackage{tabularx}
\usepackage[backend=bibtex,sorting=none,style = ieee]{biblatex}
\usepackage{lmodern}   
\usepackage{csquotes}
\usetikzlibrary{arrows.meta,calc,positioning,decorations.pathreplacing}
\pgfplotsset{compat=1.18}
\usepgfplotslibrary{groupplots}
\usetikzlibrary{external}
\tikzexternaldisable 

\makeatletter
\def\thm@space@setup{%
  \thm@preskip=2pt   
  \thm@postskip=2pt  
}
\makeatother

\newcommand*{\vm}[1]{\boldsymbol{#1}}
\newcommand*{\trans}{^{\top}}
\newcommand*{\dd}{\mathrm{d}}
\newcommand*{\nablax}{\nabla_{\vm x_i}}

\newcommand*{\nablaxj}{\nabla_{\vm x_j}}
\newcommand*{\neighs}{j \in \mathcal{N}_i}

\newcommand*{\agents}{i \in \mathcal{V}}
\newcommand*{\Ni}[1]{\vm{#1}_{\mathcal N_i}}
\newcommand*{\Nj}[1]{\vm{#1}_{\mathcal N_j}}
\newcommand*{\inds}{^{*}}

\newcommand*{\indq}{^{q}}
\newcommand*{\indjq}{^{j|q}}
\newcommand*{\indjpq}{^{j-1|q}}
\newcommand*{\st}{\operatorname{s.\!t.}}
\newcommand*{\smallspace}{\,\,}
\newcommand*{\indqn}{^{q+1}}

\newcommand*{\timetau}{\tau \in [0,\,T]}
\newcommand*{\timetausample}{\tau \in [0,\,\Delta t)}

\DeclareMathOperator*{\argmin}{\arg\, \min}

\DeclareMathOperator*{\diag}{\mathrm{diag}}

\DeclareMathOperator*{\sign}{\mathrm{sign}}
\newtheorem{theorem}{Theorem}

\newtheorem{example}{Example}

\definecolor{mycolor1}{rgb}{0.00000,0.44700,0.74100} 
\definecolor{mycolor2}{rgb}{0.85000,0.32500,0.09800} 
\definecolor{mycolor3}{rgb}{0.9290 0.6940 0.1250} 
\definecolor{mycolor4}{rgb}{0.4940 0.1840 0.5560} 
\definecolor{mycolor5}{rgb}{0.4660 0.6740 0.1880} 
\definecolor{mycolor6}{rgb}{0.3010 0.7450 0.9330} 

\pgfplotscreateplotcyclelist{WaterTankSimulation}{ 
	{mycolor1,thick},
	{mycolor1,thick},
	{mycolor2,thick,dashed},
	{mycolor2,thick,dashed},
	{black,dashed},
	{black,dashed}
}

\pgfplotscreateplotcyclelist{WaterTankExperiment}{ 
	{mycolor1,thick},
	{mycolor1,thick},
	{black,dashed},
	{black,dashed},
}

\pgfplotscreateplotcyclelist{ParameterEstimation}{ 
	{mycolor1,thick},
}

\pgfplotsset{label style={font=\footnotesize},ticklabel style={font=\footnotesize},legend style={font=\tiny}, legend cell align=left}